\documentclass[12pt,oneside]{article}
\textheight215truemm \textwidth 138truemm  \date{}
\usepackage{amsmath,amsfonts,amssymb,theorem,epsf}

 \textheight215truemm \textwidth 138truemm  \date{}
 \title{Foreword}
 \author{Alessio Corti \and Miles Reid}
 \addcontentsline{toc}{chapter}{Alessio Corti and Miles Reid: Foreword}

\newtheorem{thm}{Theorem}

{
\theorembodyfont{\rmfamily}
\newtheorem{prb}[thm]{Problem}
\newtheorem{exa}[thm]{Example}
\newtheorem{rmk}[thm]{Remark}

}

 \numberwithin{thm}{subsection}
 \numberwithin{equation}{subsection}
 \numberwithin{figure}{subsection}

 \newcommand{\rest}[1]{_{{\textstyle|}#1}} 
\newcommand{\broken}{\dasharrow}
\newcommand{\into}{\hookrightarrow}
\newcommand{\iso}{\cong}
\newcommand{\ep}{\varepsilon}
\newcommand{\fie}{\varphi}
\newcommand{\al}{\alpha}
\newcommand{\la}{\lambda}
\newcommand{\La}{\Lambda}
\newcommand{\Ga}{\Gamma}
\newcommand{\si}{\sigma}

\newcommand{\Vbar}{\overline V}
\newcommand{\pt}{\mathrm{pt}}
\DeclareMathOperator{\rk}{rank}
\DeclareMathOperator{\mult}{mult}
\DeclareMathOperator{\Pic}{Pic}
\DeclareMathOperator{\Proj}{Proj}
\DeclareMathOperator{\Exc}{Exc}
\DeclareMathOperator{\NEbar}{\overline{NE}}
\DeclareMathOperator{\Spec}{Spec}

\newcommand{\1}{^{-1}}
\newcommand{\Oh}{\mathcal{O}}
\newcommand{\sH}{\mathcal{H}}
\newcommand{\sP}{\mathcal{P}}
\newcommand{\bq}{\mathbf{q}}
\newcommand{\Q}{\mathbb{Q}}
\newcommand{\R}{\mathbb{R}}
\newcommand{\PP}{\mathbb{P}}
\newcommand{\FF}{\mathbb{F}}
\newcommand{\C}{\mathbb{C}}
\newcommand{\Z}{\mathbb{Z}}

 \makeindex
\begin{document}

 \enlargethispage{20pt}
 \setcounter{page}{0}
\maketitle
 \thispagestyle{empty}

 \markboth{\qquad Foreword \hfill}{\hfill Alessio Corti and Miles Reid
\qquad}
 \begin{abstract} This is the Foreword to the book
 \begin{quote}
Explicit birational geometry of 3-folds, edited by A.~Corti and M.~Reid,
CUP Jun 2000, ISBN: 0 521 63641 8. Papers by K. Altmann, A. Corti, A.
R. Iano-Fletcher, J.~Koll\'ar, A.~V. Pukhlikov and M. Reid.
 \end{quote}

 One of the main achievements of algebraic geometry over the last 20 years
is the work of Mori and others extending minimal models and the
Enriques--Kodaira classification to 3-folds. This book is an integrated
suite of papers centred around applications of Mori theory to birational
geometry. Four of the papers (those by Pukhlikov, Fletcher, Corti, and the
long joint paper Corti, Pukhlikov and Reid) work out in detail the theory
of birational rigidity of Fano 3-folds; these contributions work for the
first time with a representative class of Fano varieties, 3-fold
hypersurfaces in weighted projective space, and include an attractive
introductory treatment and a wealth of detailed computation of special
cases.

 \paragraph*{Contents}
 \begin{enumerate}
 \item[0.] A. Corti and M. Reid: Foreword \dotfill 1--20
 \item K. Altmann: One parameter families containing three \\ dimensional
toric Gorenstein singularities \dotfill 21--50
 \item J. Koll\'ar: Nonrational covers of $\C\PP^m\times\C\PP^n$ \dotfill
51--71
 \item A. V. Pukhlikov: Essentials of the method of maximal \\
singularities \dotfill 73--100
 \item A. R. Iano-Fletcher: Working with weighted complete \\
intersections \dotfill 101--173
 \item A. Corti, A. V. Pukhlikov and M.~Reid: Fano 3-fold \\ hypersurfaces
\dotfill 175--258
 \item A.~Corti: Singularities of linear systems and 3-fold \\ birational
geometry \dotfill 259--312
 \item M. Reid: Twenty five years of 3-folds -- an old \\ person's view
\dotfill 313--343
 \item[] {\bf Index} \dotfill 345--349
 \end{enumerate}
\end{abstract}
\clearpage

\section{Introduction}\label{sec_ex}
This volume is an integrated collection of papers working out several new
directions of research on 3-folds under the unifying theme of {\em
explicit birational geo\-metry}. Section~\ref{sec_book} summarises briefly
the contents of the individual papers.

Mori theory\index{Mori!theory|(} is a conceptual framework for studying
minimal models and the classification of varieties, and has been one of
the main areas of progress in algebraic geometry since the 1980s. It
offers new points of view and methods of attacking classical problems,
both in classification and in birational geo\-metry, and it raises many new
problem areas. While birational geo\-metry has inspired the work of many
classical and modern mathematicians, such as L.~Cremona, G.~Fano, Hilda
Hudson, Yu.~I. Manin, V.~A. Iskovskikh and many others, and while their
results undoubtedly give us much fascinating experimental material as food
for thought, we believe that it is only within Mori theory that this body
of knowledge begins to acquire a coherent shape.

At the same time as providing adequate tools for the study of 3-folds,
Mori theory enriches the classical world many times over with new examples
and constructions. We can now, for example, work and play with hundreds of
families of Fano 3-folds.\index{Fano!3-fold} {From} where we stand, we can
see clearly that the classical geometers were only scratching at the
surface, with little inkling of the gold mine awaiting discovery.

The theory of minimal models of surfaces works with nonsingular surfaces,
and the elementary step it uses is Castel\-nuovo's criterion, which allows
us to contract $-1$-curves (exceptional curves of the first kind). A chain
of such contractions leads us to a minimal surface $S$, either $\PP^2$ or
a scroll over a curve, or a surface with $K_S$ numerically
nonnegative\index{nef} (now called {\em nef}, see \ref{ssec:cone} below).
These ideas were well understood by Castel\-nuovo and Enriques a century
ago, and are so familiar that most people take them for granted. However,
their higher dimensional generalisation was a complete mystery until the
late 1970s, and may still be hard to grasp for newcomers to the field. It
involves a suitable category of mildly singular projective varieties, and
the crucial new ingredient of {\em extremal ray}\index{extremal!ray}
introduced by Mori around 1980. As we discuss later in this foreword,
extremal rays provide the elementary steps of the minimal model program
(the divisorial contractions\index{divisorial contraction} and
flips\index{flip} of the Mori category\index{Mori!category} that
generalise Castel\-nuovo's criterion) and also the definition of Mori
fibre space\index{Mori!fibre space}\index{strict Mori fibre space} (that
generalise $\PP^2$ and the scrolls), our primary object of interest.

Higher dimensional geo\-metry, like most other areas of mathematics, is
marked by creative tensions between abstract and concrete on the one hand,
general and special on the other. Contracting a $-1$-curve on a surface is
a concrete construction, whereas a Mori extremal ray and its contraction
is abstract (compare Remark~\ref{rmk:contrast}). The ``general'' tendency
in the classification of varieties, exemplified by the work of Iitaka,
Mori, Koll\'ar, Kawamata and Shokurov, includes things like
Iitaka--Kodaira dimension, cohomological methods, and the minimal model
program in substantial generality. The ``special'' tendency, exemplified
by Hudson, Fano, Iskovkikh, Manin, Pukhlikov, Mori and ourselves, includes
the study of special cases, for their own sake, and sometimes without hope
of ever achieving general status.

By {\em explicit}, we understand a study that does not rest after obtaining
abstract existence results, but that goes on to look for a more concrete
study of varieties, say in terms of equations, that can be used to bring
out their geometric properties as clearly as possible. For example, the
list of Du Val singularities by equations and Dynkin diagrams is much more
than just an abstract definition or existence result, and can be use for
all kinds of purposes. This book initiates a general program of explicit
birational geometry of 3-folds (compare Section~\ref{sec_prob}). On the
whole, our activities do not concern themselves with 3-folds in full
generality, but work under particular assumptions, for example, with
3-folds that are hypersurfaces, or have\index{quotient singularity} only
terminal quotient
singularities\index{terminal!singularities}\index{terminal!quotient sing.\
$\frac{1}{r}(1,a,r-a)$} (see below). The advantage is that we can get a
long way into current thinking on 3-folds while presupposing little in the
way of technical background in Mori theory.

Treating 3-folds and contractions between them in complete generality would
lead us of necessity into a number of curious and technically difficult
backwaters; these include many research issues of great interest to us,
but we leave them to more appropriate future publications (see however
Section~\ref{sec_prob} below). Making the abstract machinery work in
dimension $\ge4$ is another important area of current research, but the
geo\-metry of 4-folds is presumably intractable in the explicit terms that
are our main interest here.

 \section{The Mori program}\label{sec_Mori}
 This section is a gentle introduction to some of the ingredients of the
\hbox{3-fold} minimal model program, with emphasis on the aspects most
relevant to our current discussion. Surveys by Reid and Koll\'ar
\cite{ten}, \cite{kol1}, \cite{kol2} also offer introductory discussions
and different points of view on Mori theory. At a technically more
advanced level, we also recommend a number of excellent (if somewhat less
gentle) surveys: Clemens, Koll\'ar and Mori \cite{CKM}, Kawamata, Matsuda
and Matsuki \cite{KMM}, Koll\'ar and Mori \cite{KM}, Mori \cite{Mo} and
Wilson \cite{W}.

 \subsection{Terminal singularities}
 It was understood from the outset that minimal models of 3-folds
necessarily involve singular varieties (one reason why is explained in
\ref{sec_bir}). The {\em Mori category}\index{Mori!category} consists of
projective varieties with {\em terminal
singularities};\index{terminal!singularities|(} the most typical example
is\index{quotient singularity} the cyclic quotient
singularity\index{terminal!quotient sing.\ $\frac{1}{r}(1,a,r-a)$}
$\frac{1}{r}(a,r-a,1)$. Here $a$ is coprime to
$r$, and the notation means the\index{quotient singularity} quotient
$\C^3/(\Z/r)$, where the cyclic group $\Z/r$ acts by
 \[
 (x,y,z)\mapsto (\ep^ax,\ep^{r-a}y,\ep z),
 \]
 and $\ep$ is a primitive $r$th root of 1. The most common instance is
$\frac{1}{2}(1,1,1)$, the cone on the Veronese surface. The effect of
saying that this point is terminal is that if we first resolve it by
blowing up, then run a minimal model program on the resolution, we will
eventually need to contract down everything we've blown up, taking us back
to the same singularity.

There are a few other classes of terminal singularities, including isolated
hyper\-surface singularities such as $xy=f(z,t)\subset\C^4$, where
$f(z,t)=0$ is an isolated plane curve singularity, and a combination
of\index{quotient singularity} hypersurface and quotient singularity, for
example, the hyperquotient singularity obtained by dividing the
hypersurface singularity $xy=f(z^r,t)$ by the cyclic group $\Z/r$ acting by
$\frac{1}{r}(a,r-a,1,0)$. At some time you may wish to look through some
sections of Reid \cite{YPG} (especially Theorem~4.5) for a more formal
treatment. But for most purposes, the cyclic quotient singularity
\hbox{$\frac{1}{r}(a,r-a,1)$} is the main case for understanding 3-fold
geo\-metry, and if you bear this in mind, you will have little trouble
understanding this book.\index{terminal!singularities|)}

 \subsection{Theorem on the Cone}\label{ssec:cone}
\begin{figure}[ht]
 \centerline{\epsfbox{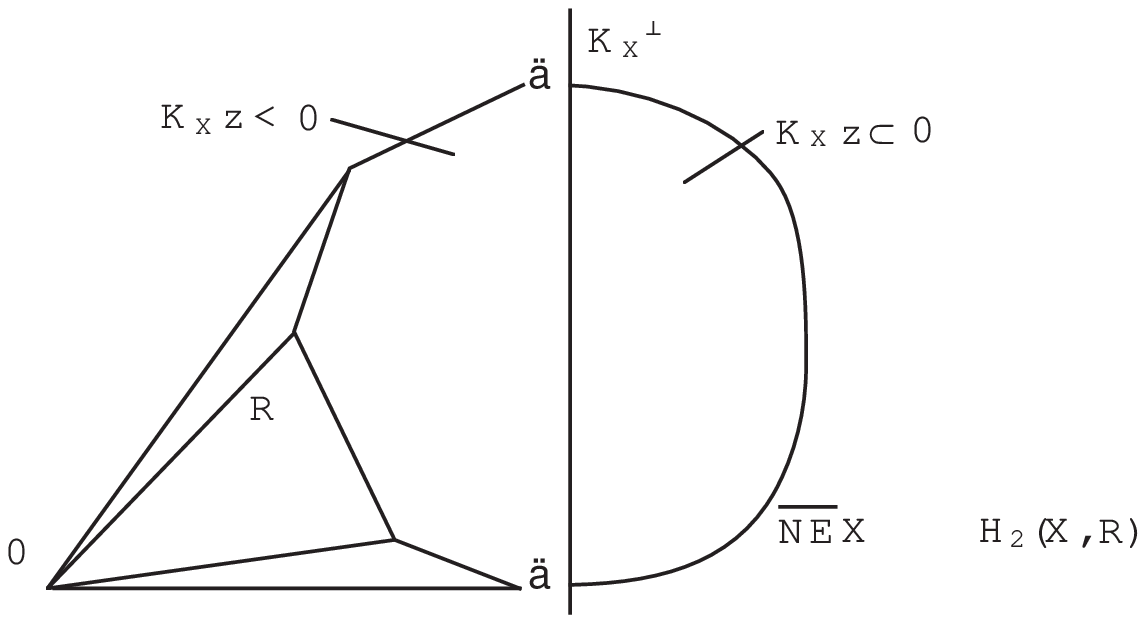}}
 \caption{The Mori cone: $\NEbar X$ is locally rational polyhedral in
$K_Xz<0$}
 \label{fig:NEbar}
 \end{figure}

 \index{Mori!cone|(}
 The {\em Mori cone\/} $\NEbar X$ (see Figure~\ref{fig:NEbar}) is
probably the most profound and revolutionary of Mori's contributions to
3-folds.
An $n$-dimensional projective variety $X$ over $\C$ is a $2n$-dimensional
oriented compact topological space, and its second homology group
$H_2(X,\R)$ is a finite dimensional real vector space. Every algebraic
curve $C\subset X$ can be triangulated and viewed as an oriented 2-cycle,
and thus has a homology class $[C]\in H_2(X,\R)$. Then by definition
$\NEbar X$ is the closed convex cone in $H_2(X,\R)$ generated by the
classes $[C]$ of algebraic curves $C\subset X$. You can think of this as
follows: $H_2(X,\R)$ is a property of the topological space $X$, whereas
the structure of $X$ as a projective algebraic variety provides the extra
information of the Mori cone $\NEbar X\subset H_2(X,\R)$.

The shape of $\NEbar X$ contains information about linear systems and
embeddings $X\into\PP^N$. Taking intersection number $D\cdot C$ with a
divisor, or evaluating $\al\cap[C]$ with a cohomology class $\al\in
H^2(X,\R)$ (say, the first Chern class of a line bundle $L$) defines a
linear form on $H_2(X,\R)$. We say that $D$ or $\al$ is {\em
nef}\index{nef} if this linear form is $\ge0$ on $\NEbar X$; that is, a
divisor $D$ is nef if $D\cdot C\ge0$ for every curve $C\subset X$. Under
an embedding, every algebraic curve must have positive degree; it is known
that, under rather mild assumptions, $X$ is projective if and only if
$\NEbar X$ is a genuine cone with a point.

To state Mori's theorem, we assume that the canonical divisor class $K_X$,
(or equivalently, the first Chern class of the cotangent bundle) makes
sense as a linear form on $H_2(X,\R)$. This is a mild extra assumption on
$X$, that certainly holds if $X$ is nonsingular or has at worst quotient
singularities.\index{quotient singularity} The theorem on the cone then
says that $\NEbar X$ is a rational polyhedral cone in the half-space of
$H_2(X,\R)$ on which $K_X$ is negative. This theorem is particularly
powerful for Fano varieties, defined by the condition that $-K_X$ is
ample: for these, the entire cone $\NEbar X$ is contained in $K_Xz<0$, so
that $\NEbar X$ is a finite rational polyhedral cone.

 \subsection{Extremal rays and the contraction theorem}\label{sec_con}
 \index{extremal!ray|(}
 Mori theory applies mainly to varieties with $K_X$ not nef. This
condition says that $K_XC<0$ for some curve $C$, or that the part of the
cone $\NEbar X$ in the half-space $K_Xz<0$ is nonempty. Since this part of
the cone is locally rational polyhedral, it follows that, if $K_X$ is not
nef,\index{nef} $\NEbar X$ has at least one extremal ray $R$ with
$K_X\cdot R<0$. Here an {\em extremal ray} is just a half-line $R=\R_+
z\subset\NEbar X$ that is extremal in the sense of convex geo\-metry (that
is,
 \[
 \text{$z_1,z_2\in\NEbar X$ and $z_1+z_2\in R$}
 \ \implies\ z_1,z_2\in R).
 \]

 Let $R\subset\NEbar X$ be an extremal ray with $K_X\cdot R<0$. Then there
exists a contraction morphism
 \[
 f_R\colon X\to Y,
 \]
characterised by the property that a curve $C\subset X$ is mapped to a
point if and only if $C\in R$ (more precisely, the class of $C$). The
morphism $f_R\colon X\to Y$ is called a {\em Mori contraction} or an {\em
extremal contraction}. It is determined by the extremal ray $R$, and has
categorical properties such as $-K_X$ relatively ample and $\rho(X/Y)=1$
that turn out to be surprisingly strong: for example $-K_X$ ample puts us
in a position where vanishing results based on Kodaira vanishing kill
almost all the cohomology.

 The cone and contraction theorems are proved in Koll\'ar and Mori
\cite{KM}; on the whole, we can get by without reference to the
technicalities of the proof, and you may prefer to take these results on
trust for now.\index{Mori!cone|)}

 \subsection{Types of extremal rays}\label{sec_cl}
 The next step is the case division on the dimension of the image $Y$ and
of the exceptional locus of the contraction morphism $f_R\colon X\to Y$,
called the {\em classification of extremal rays} (or {\em rough
classification}). The cases when the contraction $f_R\colon X\to Y$ has
$\dim Y<\dim X$ lead to the definition of\index{Mori!fibre space} Mori
fibre space and Fano varieties\index{Fano!variety} discussed in
\ref{sec_Mfs}. In the other cases, we are dealing with birational
modifications of $X$, and, as we see in \ref{sec_bir}, the aim is to
proceed inductively towards a minimal model, as in the classical case of
surfaces.

\begin{rmk}\label{rmk:contrast} Note the contrast with the classical case:
for surfaces, the thing we contract is a geo\-metric locus. We find a
$-1$-curve $C$ and establish that it can be contracted in terms of a
neighbourhood of $C$. In contrast, Mori theory in dimension $\ge3$ works
primarily in terms of {\em categorical definitions} and {\em existence
theorems}: the thing to be contracted is an extremal ray $R$ of $\NEbar X$
(the definition of which uses the totality of curves on $X$). The proof of
the general theorems saying that $R$ is contractible by a morphism $f_R$
makes sophisticated use of numerical and cohomology vanishing properties
of $X$.

The geometric nature of the contraction is only studied as a second step;
even basic things such as the geometric locus that is contracted or even
the dimension of the image cannot be anticipated. $f_R$ may be birational,
a proper fibre space, or the constant morphism to a point. This curious
inversion of thinking is another of Mori's characteristic contributions to
the subject, and the logic still comes as a surprise to anyone knowing a
traditional treatment of the classification of surfaces. After all,
Castelnuovo and Enriques could scarcely have guessed that (i)~contracting
a $-1$-curve, (ii) projecting a geo\-metrically ruled surface to its base
curve, and (iii) the constant map of $\PP^2$ to a point would find a
unified treatment as extremal contractions,\index{extremal!contraction}
and that this idea, however outlandish it might appear at first sight,
would lay the foundations of all future work in
classification.\index{extremal!ray|)}
 \end{rmk}

 \subsection{Birational modifications: divisorial contractions, flips and
the minimal model program} \label{sec_bir} \index{divisorial contraction}
\index{flip}
 The extremal contractions\index{extremal!contraction} that are most
similar to contracting a $-1$-curve on a non\-singular surface
(Castelnuovo's criterion) are the {\em divisorial contractions}. Here the
case assumption is that $f_R\colon X\to Y$ is birational, and contracts a
divisor of $X$ to a locus of $Y$ of co\-dimension $\ge2$. The categorical
properties of $f_R$ then guarantee automatically that the exceptional
locus of $f_R$ is an {\em irreducible} divisor, and that $Y$ has terminal
singularities.\index{terminal!singularities} This is the point at which
terminal singularities force themselves on our attention: even if $X$ is
nonsingular, $Y$ may be singular. Because $Y$ is still in the Mori
category,\index{Mori!category} we can repeat the same game starting from
$Y$.

 \index{small!contraction} The other birational case is when $f_R$ is {\em
small}, that is, every component of $\Exc f_R$ has codimension $\ge2$; in
this case there cannot be any cohomology class in $H^2(X,\R)$ that
corresponds to the canonical divisor of $Y$, so that $Y$ can {\em never}
have terminal singularities.\index{terminal!singularities} (If such a
class existed, its pullback to $X$ would coincide with $K_X$, which would
then be numerically trivial on the fibres of $f_R$. This contradicts
$-K_X$ ample, the defining property of a Mori extremal contraction.)

 Because $Y$ is no longer in the Mori category, the minimal model program
cannot just continue inductively from $Y$. The subject was stuck at this
point for a few years in the 1980s, before Mori proved the 3-fold flip
theorem:\index{flip!theorem} there is a {\em flip}\index{flip}
 \begin{equation}
 \renewcommand{\arraycolsep}{0.1em}
 \renewcommand{\arraystretch}{1.15}
 \begin{array}{ccccc}
 X \kern-0.12em &&\stackrel{t_R}{\broken}&& \kern-0.28em X^+ \\
 & \searrow && \swarrow \\
 && Y
 \end{array}
 \label{eq_flip}
 \end{equation}
where $X^+\to Y$ is another birational map from a 3-fold $X^+$,
characterised by the property that $K_{X^+}$ is ample over $Y$. In other
words, the birational map $t_R\colon X\broken X^+$ cuts out from $X$ a
finite number of curves on which $K_X$ is negative, and in their place
glues back into $X^+$ a finite number of curves on which $K_{X^+}$ is
positive. The definition of flip may seem somewhat obscure, but many nice
attributes of $X^+$ follow from it; in particular, the morphism $X^+\to Y$
is also small,\index{small!contraction} and $X^+$ again has terminal
singularities,\index{terminal!singularities} so is in the Mori
category.\index{Mori!category} In dimension $\ge4$, the existence of the
flip diagram (\ref{eq_flip}) is called the {\em flip
conjecture};\index{flip!conjecture} this seems to be one of the most
intractable problems in the subject.

Divisorial contractions and Mori flips\index{flip} are the elementary
steps in the Mori minimal model program. A sequence of these leads after a
finite number of steps to a variety $X'$, which is either a {\em minimal
model}, that is, a variety with $K_{X'}$ nef,\index{nef} or a Mori fibre
space\index{Mori!fibre space} $f\colon X'\to S$. 

 \subsection{The definition of Mori fibre space} \index{Mori!fibre space}
 \index{strict Mori fibre space} \label{sec_Mfs}
 We now discuss the remaining cases in the classification of extremal
rays,\index{extremal!ray} when the contraction $f_R\colon X\to Y$
maps to a smaller dimensional variety, that is, $\dim Y<\dim X$. Then
$f_R$ (or $X$ itself) is called a {\em Mori fibre space} (Mfs). Note that,
following Iitaka and Ueno, we say {\em fibre space} to mean a morphism
$f\colon X\to Y$, often assumed to have connected fibres and $Y$ normal,
possibly with varying fibres, singular fibres, even fibres of different
dimensions; this is not to be confused with the much stricter notion of
fibre bundle.

The cases when $Y$ is a surface and $X\to Y$ is a conic bundle\index{conic
bundle} (that is, the general fibre is a conic) or when $Y$ is a curve and
$X\to Y$ a fibre space of del Pezzo surfaces are the natural analogues of
ruled surfaces. For the logical framework of Mori theory, we include in
the definition of Mori fibre space the case that the contraction
$f_R\colon X\to Y=\pt.$ is the constant map to a point: then the morphism
$f_R$ is trivial, but its categorical properties include the fact that
$-K_X$ is ample, and $\Pic X$ has rank~1. In this case $X$ is called a
{\em Fano $3$-fold\/};\index{Fano!3-fold} in contrast to the classical
terminology, we allow $X$ to be singular.

\subsection{Biregular geometry versus birational geometry}\label{ssec:bir}

 The dividing line between biregular and birational geometry has changed
through the generations, and is possibly still open to debate. The Italian
school worked primarily in birational terms, and Zariski and Weil used
birational ideas (at least in part) in setting up foundations for
biregular geometry. The modern view, with scheme theory firmly established
as the foundation, constructs birational geometry within this biregular
framework. Thus, while the dichotomy between surfaces having nonvanishing
plurigenera and ruled surfaces (or ``adjunction terminates'') is manifestly
birational, we no longer think of it as the primary result of
classification, but derive it from biregular results. This new view was
instrumental in the success of Mori theory.

 When we run a Mori minimal model program on a given 3-fold $V$, the end
product is either a minimal model $X$ with $K_X$ nef,\index{nef} or a Mori
fibre space,\index{strict Mori fibre space} typically, a Fano 3-fold $X$
or a conic bundle\index{conic bundle} over a surface $X\to S$. The
properties that define the 3-fold $X$ are biregular in nature, so that we
view $X$ as a biregular construction. {From} this point of view, the proof
of classification should also be considered a biregular activity, since the
point is to prove that a given minimal 3-fold $X$ has the right
plurigenera and Kodaira dimension. Our conclusion is that birational
geo\-metry begins with the question of birational maps between different
Mori fibre spaces.\index{Mori!fibre space}\index{Mori!theory|)}

\section{What this book contains} \label{sec_book}

This section discusses briefly the papers in this book, and their
contribution to the above program of study. The papers are:
 \begin{enumerate}
 \renewcommand{\labelenumi}{(\arabic{enumi})}
 \item K. Altmann: One-parameter families containing
three-dimensional\newline toric Gorenstein
singularities\index{toric!Gorenstein singularities}
 \item J. Koll\'ar: Nonrational covers of $\PP^m \times \PP^m$
 \item A. V. Pukhlikov: Essentials of the method of maximal singularities
 \item A. R. Iano-Fletcher: Working with weighted complete intersections
 \item A. Corti, A. Pukhlikov and M. Reid: Fano 3-fold hypersurfaces
 \item A. Corti: Singularities of linear systems and 3-fold birational
geo\-metry
 \item M. Reid: Twenty five years of 3-folds, an old person's view
 \end{enumerate}

 Klaus Altmann's paper~(1) is a study of the\index{deformation theory}
deformation theory of\index{toric!Gorenstein singularities} toric
Gorenstein 3-fold singularities. It relates to the classification of
3-fold flips\index{flip} as follows: we know that any Mori flip diagram
(\ref{eq_flip}) can be obtained from a $\C^\times$ action on a 4-fold
Gorenstein singularity $0\in A$ by taking the quotient by the $\C^\times$
action in different interpretations -- the so-called {\em variation of
geo\-metric invariant theory quotient}, see Dolgachev and Hu~\cite{DH},
Reid~\cite{R2} and Section~\ref{sec_flips} below. Moreover, the general
anticanonical\index{anticanonical!divisor} divisor $S\in|{-}K_X|$ (the {\em
general elephant\/}) is a surface with only Du Val singularities,
according to Koll\'ar and Mori \cite{KM1}, Theorem~1.7. Its inverse image
in $A$ is a $\C^\times$ cover $B\to S$, and is a hyperplane section
$B\subset A$, so that $A$ can be viewed as a 1-parameter deformation of
$B$. It frequently happens that $S$ is of type $A_n$, and then $B$ is
toric,\index{toric!variety} so that Altmann's theory applies in many cases
to give a classification of 3-fold flips. Altmann's previous work
\cite{Al} used the notion\index{Minkowski!sum} of Minkowski decomposition
of polytopes to give a complete treatment of the deformation of isolated
3-fold toric Gorenstein singularities;\index{toric!Gorenstein
singularities} in the present paper, he shows how to modify his method to
the case of toric varieties having singularities in codimension 2.

J\'anos Koll\'ar's paper~(2) provides a new method of proving
irrationality, adding to the known collection of rationally
connected\index{rationally connected} varieties that\index{nonrational
variety} are not rational: finite covers of $\PP^m\times\PP^n$ with
ramification divisor of large enough degree in one factor, and
hypersurfaces in $\PP^m\times\PP^n$ of large enough degree. His technique
involves reduction to characteristic $p$, and a rather clever and
surprising analysis of the stability of the tangent bundle in
characteristic $p$. In fact, he proves the slightly more general
structural property that these varieties are not even ruled. In the case
of conic bundles,\index{conic bundle} these results are spectacularly
close to the conjectural\index{rationality problem} bound for rationality
(compare, for example, paper~(2), Remark~1.2.1.1 with Corti's paper~(6),
4.10 and 4.11). This provides the strongest confirmation to date of the
conjectures on conic bundles,\index{conic bundle} in a numerical range
that is inaccessible to all other methods.

 \index{birational rigidity}
 The papers~(3)--(6) form a connected suite of papers around the subject of
{\em birational rigidity}. The notion, discussed in more detail in
Section~\ref{sec_birig} below, originates in the famous result of
Iskovskikh and Manin \cite{IM} that a non\-singular quartic 3-fold
$X_4\subset\PP^4$ has no birational maps to Fano varieties (other than
isomorphisms to itself). Pukhlikov's paper~(3) describes his important
simplification and elaboration of Iskovskikh and Manin's treatment. This
paper is partly based on notes of lectures given at the 1995--96 Warwick
algebraic geo\-metry symposium and the preprint, with its clear treatment
of the Russian methods, strongly stimulated our collaboration in the joint
paper~(5). The different approach in Pukhlikov's papers also offers a
useful ideological and practical counterweight to the methods of Corti's
paper~(6).

Our long joint paper~(5) is the real heart of this book. In it, we carry
out a substantial portion of a program of research on birational rigidity,
treating the {\em famous 95\/} families\index{famous@``famous 95''} of Fano
3-fold\index{Fano!3-fold} weighted hypersurfaces. We refer to
Section~\ref{sec_birig} and the introduction to paper~(5) for further
discussion of birational rigidity.

Anthony Iano-Fletcher's paper~(4) is a well written tutorial introduction
to weighted projective spaces and their subvarieties. This paper has been
available for many years as a Max Planck Institute preprint, and is widely
quoted in the literature; it contains many very useful results and methods
of calculation, including one derivation of the list of the famous 95
hyper\-surfaces, and thus forms an essential prerequisite for paper~(5).

Corti's paper~(6) contains a detailed introduction to the Sarkisov
program.\index{Sarkisov!program} It develops and applies powerful new
methods to quantify and analyse the singularities of linear systems,
clarifying and providing technical alternatives to the methods initiated
by Iskovskikh and Manin based on the study of the resolution graph. The
new ideas are based on the Shokurov connectedness\index{Shokurov
connectedness} principle in log birational geo\-metry,\index{log!category}
and seem to provide the most powerful currently known technique to exclude
birational maps between Mori fibre spaces.\index{strict Mori fibre space}
The results are applied to give rigidity\index{birational rigidity}
criteria for Mori fibre spaces in a number of cases, and our joint
paper~(5)\index{log!surface method} also appeals to them for one or two
technical points.

Reid's historical paper~(7) is a {\em Heldenleben} that needs no
introduction.

 \section{Mori theory and birational geometry} \label{sec_Sark}
 \index{Mori!theory}
 Following our introductory remarks on Mori theory in
Section~\ref{sec_Mori}, we now give a brief introduction to our view of
birational geometry, including the Sarkisov
program\index{Sarkisov!program} and birational
rigidity.\index{birational rigidity}

 \subsection{Fano-style projections}\index{projection}
 Fano based his treatment of the 3-folds $V_{2g-2}\subset\PP^{g+1}$ for
$g\ge7$ on the idea of constructing a birational map by projection from a
suitably chosen centre. Typically, the double projection of $V$ from a
line $L$ involves a diagram
 \begin{equation}
 \renewcommand{\arraystretch}{1.3}
 \renewcommand{\arraycolsep}{0.25em}
 \begin{matrix}
 V' & \broken & V'' \\
 \downarrow && \downarrow \\
 V && W,
 \end{matrix}
 \label{eq_Fano}
 \end{equation}
where $V'\to V$ is the blowup of $L\subset V$, the map $V'\broken V''$
flops the lines meeting $L$ (in good cases, finitely many lines with
normal bundle of type $(-1,-1)$), and $V''\to W$ contracts the surface
$E\subset V''$ swept out by conics meeting $L$ to a curve $\Ga\subset W$.
Fano thought of the map $V\broken W$ as the rational map defined by linear
projection, and factoring it in biregular terms was not his primary
concern. 

For us, on the other hand, it is important to view Fano's projection as a
general construction in the Mori category:\index{Mori!category} $V'\to V$
is an extremal extraction,\index{extremal!extraction} $V'\broken V''$ a
rational map that is an isomorphism in codimension~1 (in good cases, a
composite of classical flops), and $V''\to W$ the contraction of an
extremal ray.\index{extremal!ray} All 4 of the varieties in
(\ref{eq_Fano}) are in the Mori category,\index{Mori!category} and the two
morphisms are contractions of extremal rays.\index{extremal!ray}

 \begin{rmk}
We take the opportunity to clear up a possible source of confusion that
occurs throughout the subject: in Fano's case, the {\em single} projection
from $L$ contracts the flopping lines by a morphism $V'\to\Vbar$ to a
variety $\Vbar$ having (in good cases) only 3-fold ordinary double points;
we think of $\Vbar$ as the {\em midpoint} of the construction of the link
$V\broken W$. It is a Fano variety\index{Fano!3-fold} in some sense, since
it has terminal singularities\index{terminal!singularities} and
$-K_{\Vbar}$ is ample, but it is {\em not in the Mori category}, because
it is not $\Q$-factorial:\index{factorial@$\Q$-factorial} the exceptional
scroll over $L$ maps to a divisor in $\Vbar$ that is not Cartier and not
$\Q$-Cartier at the nodes.
 \end{rmk}

 \subsection{Sarkisov links}
 Sarkisov links\index{link!of Type~I--IV} play the role of
``elementary transformations'' for birational maps between Mori fibre
spaces.\index{Mori!fibre space}\index{strict Mori fibre space} The {\em
Sarkisov program}\index{Sarkisov!program} factors an arbitrary birational
map
$X\broken Y$ between the total spaces of Mori fibre spaces $X\to S$ and
$Y\to T$ as a chain of such links (see Corti's paper~(6) and \cite{Co} for
details).

A general {\em Sarkisov link}\index{link} is given by a diagram in
the Mori category\index{Mori!category} that is a variation on
(\ref{eq_Fano}), but with general extremal
contractions\index{extremal!contraction} allowed as the morphisms. As
discussed in much more detail in Corti \cite{Co} and paper~(6), if we
start from a Mori fibre space\index{Mori!fibre space}
$X\to S$, any link $(X/S)\broken(Y/T)$ is given by one of the following
constructions. First, we replace $X\to S$ by a new morphism $X_1\to S_1$
having $\rk N^1(X_1/S_1)=2$: for this, either
 \begin{enumerate}
 \renewcommand{\labelenumi}{(\roman{enumi})}
 \item blow up $X$ by an extremal blowup\index{extremal!extraction}
$X_1\to X$, and let $X_1\to S_1$ be the composite $X_1\to X\to S=S_1$; or
 \item contract the base by an extremal
contraction\index{extremal!contraction} $S\to S_1$, and let
$X_1\to S_1$ be the composite $X_1=X\to S\to S_1$.
 \end{enumerate}

In either case $N^1(X_1/S_1)=\R^2$, and $\NEbar(X_1/S_1)$ has an initial
extremal ray\index{extremal!ray} corresponding to the given morphism
$X_1\to X$ or $X_1\to S$. This sets up a restricted type of minimal model
program called a {\em $2$-ray game}:\index{2-ray game} because a cone in
$\R^2$ is just a ``wedge'', it has a far side that is a (pseudo-) extremal
ray (possibly not of Mori type). The case that leads to a link is when the
minimal model program runs to completion in the Mori
category:\index{Mori!theory} the far ray can be contracted, and possibly
after a chain of inverse flips,\index{inverse flip} flops and flips, the
minimal model program ends with a divisorial contraction\index{divisorial
contraction} $X^{(n)}\to Y\to S_1=T$ or a contraction of fibre type
$X^{(n)}=Y\to T\to S_1$. There are two possible ways of starting the
construction, and two ways of ending it, leading to Sarkisov links of
Type~I--IV.\index{link!of Type~I--IV}

Existence and uniqueness: the 2-ray game\index{2-ray game} is entirely
determined by the initial step $X_1/S_1$. It may happen that the initial
step $X_1/S_1$ does not construct a link\index{bad link} -- either because
an inverse flip demanded by the 2-ray game\index{2-ray game} does not
exist or has worse than terminal singularities,\index{terminal!singularities} or because the final divisorial\index{divisorial
contraction} or fibre type contraction falls out of the Mori
category.\index{Mori!category} See paper~(5), Section~5.5 and~7.6 for
examples.

For surfaces over an algebraically closed field, links are the following
familiar transformations: the blowup taking $\PP^2$ to $\FF_1$, its
inverse contraction $\FF_1\to\PP^2$, the well known elementary
transformations $\FF_k\broken\FF_{k\pm1}$ between scrolls over $\PP^1$,
and the ``exchange of factors'' of $\FF_0=\PP^1\times\PP^1$ (in other
words, the identity on $\FF_0$, but viewed as exchanging its two
projections).\index{projection} These are exactly the elementary steps in
Castelnuovo's proof of Max Noether's theorem (discussed in \ref{ssec_Noe}
below).

 \subsection{The Sarkisov program}\index{Sarkisov!program}

 We explain in conceptual terms how to factor (or
``untwist'')\index{untwisting|(} a birational map $\fie\colon X\broken X'$
\hbox{between} Mori fibre spaces\index{Mori!fibre space} $X\to S$ and
$X'\to S'$ as a chain of Sarkisov links.\index{link} Untwisting
is a constructive descending induction: let $\sH'$ be a very ample
complete linear system on $X'$, chosen at the outset and kept fixed
throughout. Following the classical ideas of Cremona, Noether and Hudson,
consider the linear system\index{mobile linear system} $\sH=\fie\1_*\sH'$
on $X$ obtained as the birational transform of $\sH'$. Untwisting is the
story of how we reduce the singularities and the degree of $\sH$ to make
$\fie$ an isomorphism.

 First we prove the {\em Noether--Fano--Iskovskikh
inequalities},\index{NFI@{Noether--Fano--Iskovskikh inequalities}} that
serve as a sensitive detector to locate the initial step $X_1/S_1$ of a
Sarkisov link\index{link} if $\fie$ is not already an
isomorphism: this is either a blowup of a maximal
singularity\index{maximal!singularity} of $\sH$, or a way of viewing $X$
over a different base to make $\sH$ look simpler. Next, when needed to
factor a given map, the Sarkisov link $\psi\colon(X/S)\to(Y/T)$ with given
initial step $X_1/S_1$ always exists.

At the start of the proof, we set up a discrete invariant of $\fie$, its
{\em Sarkisov degree}\index{Sarkisov!degree} $\deg\fie$. We only explain
this for a Fano 3-fold\index{Fano!3-fold} $X$, when $-K_X$ is a $\Q$-basis
of $\Pic X$. Then we set $\deg\fie=n$, where $n$ is the positive rational
number for which $\sH\subset|{-}nK_X|$. We prove that the Sarkisov
link\index{link} $\psi\colon X\broken Y$ provided by the NFI
inequalities decreases the Sarkisov degree, in the sense that the
composite map $\fie\psi\1\colon Y\broken X'$ between $Y/T$ and $X'/S'$ has
 \[
 \deg\fie\psi\1<\deg \fie.
 \]
We say that $\fie\psi\1$ is an {\em untwisting}\index{untwisting} of
$\fie$ by $\psi$. The factorisation theorem then follows by descending
induction on the Sarkisov degree.\index{Sarkisov!degree} Of course, we are
glossing over many subtle points, including the definition of Sarkisov
degree for a strict Mori fibre space\index{Mori!fibre space} $X\to S$, the
verification that untwisting\index{untwisting} by a link decreases
$\deg\fie$, and that a chain of untwistings must terminate. For the
details, see paper~(6).

\subsection{A classical example}\label{ssec_Noe}
We content ourselves with illustrating how these ideas work in the most
famous case of all, a birational map from $X=\PP^2$ to $X'=\PP^2$. Max
Noether's inequality states that there are 3 points $P_1,P_2,P_3$ of
$\PP^2$ (possibly infinitely near), such that
 \begin{equation}
 m_1+m_2+m_3\ >\ n=\deg\fie, \quad\text{where $m_i=\mult_{P_i}\sH$.}
 \label{eq_Noether}
 \end{equation}
Here $\mult_{P_i}\sH$ means the multiplicity of a general element of the
linear system\index{mobile linear system} $\sH$ at $P_i$. In the general
case, when $P_1,P_2,P_3$ are distinct noncollinear points, we can choose
coordinates so that these are the three coordinate points $(1,0,0)$,
$(0,1,0)$, $(0,0,1)$, and it is easy to check that untwisting by the
standard quadratic Cremona involution\index{involution!Cremona}
 \[
 \psi\colon (x_0:x_1:x_2)\broken (x_1x_2:x_0x_2:x_1x_2)
 \]
decreases the degree. Indeed, $\psi\1(\mathrm{line})$ is a conic passing
through the $P_i$, so that after untwisting, the degree becomes
 \[
 2n-\sum m_i<n.
 \]

The gap in this argument is that the points $P_i$ can be infinitely near.
This can happen in two or three different ways, and in most of these
cases, we can still construct a quadratic transformation centred on a
suitably chosen coordinate triangle to untwist our map and decreases its
degree. However, there are cases in which no single quadratic Cremona
transformation decreases the degree.

This is the starting point of Castelnuovo's proof of Noether's theorem, a
direct precursor of the Sarkisov program.\index{Sarkisov!program} Whatever
infinitely near points there may be, there always exists one point
$P\in\PP^2$ with $m_P>\frac n3$. Of course, this follows from
(\ref{eq_Noether}), but it is also easy to prove directly by an easy
argument in the spirit of ``termination of adjunction'' (see paper~(6),
Theorem~2.4, where the inequality and the argument to prove it are
generalised to any Mori fibre space).\index{Mori!fibre space} Blowing up
this point by $\FF_1\to\PP^1$ is the first link in a Sarkisov chain. It
untwists because the Sarkisov degree\index{Sarkisov!degree} measures
divisors on $\PP^2$ in terms of $-K_{\PP^2}=\Oh(3)$, but measures relative
divisors on $\FF_1$ in terms of the relative $-K_{\FF_1/\PP^1}=\Oh(2)$
(\cite{Co}, 1.3 contains a more detailed description).
\index{untwisting|)}

\subsection{Birational rigidity}\label{sec_birig}
 \index{birational rigidity|(}
 A Fano variety\index{Fano!variety} $X$ is {\em birationally rigid\/} if
for every Mori fibre space\index{Mori!fibre space} $Y\to S$, the existence
of a birational equivalence $X\broken Y$ implies that $Y\iso X$. Once the
Sarkisov program\index{Sarkisov!program} is established (that is, as yet
only in dimension $\le3$), it is equivalent to say that there are either
no Sarkisov links\index{link} out of $X$, or only self-links
$X\broken X$.

The notion of birational rigidity for strict Mori fibre spaces is also
well studied. However, the definition (paper~(6), Definitions~1.2--3) is
subtle and somewhat confusing, largely because it treats Mori fibre
spaces\index{Mori!fibre space} up to {\em square birational
equivalence}.\index{square birational map} It covers important results of
Sarkisov\index{Sarkisov!theorem on conic bundles} on conic
bundles\index{conic bundle} and Pukhlikov on del Pezzo fibre spaces.

The main result of paper~(5) states that a general member of any of the 95
families\index{famous@``famous 95''} of Fano 3-fold\index{Fano!3-fold}
hypersurfaces is rigid. We show that any birational map $\fie\colon
X\broken Y$ factors as a chain of Sarkisov links\index{link}
$X\broken X$, followed by an isomorphism $X\iso Y$. The paper is written
to be essentially selfcontained. The Sarkisov
program\index{Sarkisov!program} is introduced in a context where it is
made simpler by various special circumstances. However, at one point, we
rely on a technical statement (Theorem~5.3.3) that is proved in Corti's
paper~(6) (although it goes back essentially to Iskovskikh and Manin, and
can be proved by the technique of Pukhlikov's paper~(3)). A substantial
part of paper~(5) is devoted to the classification of links $X\broken X$,
which we also discuss in Section~\ref{sec_links} below. It was rather
surprising for us to discover that all the links can be described
explicitly in terms of just two basic constructions in commutative
algebra, which are natural generalisations of the classical Geiser and 
Bertini involutions of cubic surfaces.\index{involution!Geiser and 
Bertini}

\subsection{Beyond rigidity}
\index{rigid!boundary|(}
We now know a handful of varieties having precisely two models as Mori
fibre space,\index{Mori!fibre space} either two Fano 3-folds (Corti and
Mella \cite{CM}), or two del Pezzo fibrations\index{del Pezzo!fibration}
(Grinenko \cite{Gr}), or one of each (in the case $X_{3,3}\subset
\PP(1,1,1,1,2,2)$ suggested by Grinenko and Pukhlikov). These varieties
are not actually rigid, but nearly so. These and other examples suggest
the following idea. Say that a birational map $\fie\colon X\broken X'$
between Mori fibre spaces $X\to S$ and $X'\to S'$ is a {\em square
equivalence}\index{square birational map} if the following two conditions
hold:
 \begin{enumerate}
 \renewcommand{\labelenumi}{(\arabic{enumi})}
 \item There is a birational map $S\broken S'$ making the obvious diagram
commute. This condition is equivalent to saying that the generic fibre of
$X\to S$ is birational to the generic fibre of $X'\to S'$ under $\fie$.

 \item The birational map of (1), from the generic fibre of $X\to S$ to a
generic fibre of $X'\to S'$, is in fact biregular.
 \end{enumerate}
We define the {\em pliability}\index{pliability} of a Mori fibre space
$X\to S$ as the set
 \[
 \sP(X/S)=\bigl\{\text{Mfs}\; Y\to T \bigm| X \;\text{is birational
to}\;Y\bigr\} /\text{square equivalence}.
 \]

Our general philosophy can be described as follows. We would like to
describe the pliability $\sP(X/S)$ of a Mori fibre space\index{Mori!fibre
space} $X\to S$ in terms of its biregular geometry. There are reasons for
thinking that $\sP(X/S)$ will often have a reasonable description, say as
a finite set or a finite union of algebraic varieties. A case division
based on the various possibilities for the size of $\sP(X/S)$ can be used
as a further birational classification of Mori fibre spaces. Note that our
general philosophy can only be stated in the language of Mori
theory.\index{Mori!theory}

Although the search for counterexamples to the L\"uroth theorem played an
important part in kick-starting the study of 3-folds around 1970, the last
25 years have seen little progress on criteria\index{rationality problem}
for rationality and unirationality, and these questions seem likely to
remain intractable in the foreseeable future. The wealth of new examples
in Mori theory must in any case cast doubt on the position of these
problems as central issues in birational geometry; they date back after
all to the golden days of innocence, when the classics (from Cremona
through to Iskovskikh) had never really met the typical examples of
birational geometry.\index{rationally connected} Rationally connected
seems to be the most useful modern replacement for (uni-)rationality,
since it is robust on taking surjective image or under deformations, and
there are good criteria for it. It is not that we are hostile to the
rationality problem; rather, since we are committed to the classification
of 3-folds, we need a theoretical framework capable of accomodating all
varieties, with all their wealth of individual behaviour. Our suggested
notion of pliability\index{pliability} is a tentative step in this
direction.\index{birational rigidity|)}\index{rigid!boundary|)}

 \section{Some open problems} \label{sec_prob}

 \subsection{Explicit birational geometry}\label{sec:explicit} In a sense,
we can define explicit birational geometry as the concrete study (including
classification) of
 \begin{enumerate}
 \renewcommand{\labelenumi}{(\arabic{enumi})}
 \item divisorial contractions,\index{divisorial contraction}
flips,\index{flip} and Mori fibre spaces;\index{Mori!fibre space}
 \item the way divisorial contractions and flips combine to form the links
of the Sarkisov program.\index{Sarkisov!program}
\end{enumerate}

The model is \cite{YPG}, Theorem~4.5, which classifies all 3-fold terminal
singularities\index{terminal!singularities} as a reasonably concrete,
explicit and finite list of families. It seems reasonable to hope that
Mori flips,\index{flip} divisorial contractions,\index{divisorial
contraction} Fano 3-folds\index{Fano!3-fold} and Sarkisov
links\index{link} between 3-fold Mori fibre spaces will
eventually succumb to a similar treatment. The overall aim is to make
everything else tractable in the same sense as the terminal singularities.
This program can be expected to provide useful employment for algebraic
geometers over several decades.

\subsection{Divisorial contractions}

 \begin{prb} Fix a 3-fold terminal singularity\index{terminal!singularities} $P\in Y$ (an analytic germ). Write down all 3-fold
divisorial contractions\index{divisorial contraction} $f\colon(E\subset
X)\to(P\in Y)$ by explicit equations.
 \end{prb}

This seems a rather difficult problem in general. To do something useful
with it, it is important to realise that, whereas you are free to pick
your favourite singularity $P\in Y$, it is then your responsibility to
classify all possible extremal blowups\index{extremal!extraction} $X\to
Y$ of $P$. The few known cases of this problem are $P\in
Y=\frac{1}{r}(1,-1,a)$ treated\index{Kawamata blowup} by Kawamata
\cite{Ka}, the ordinary node (see paper~(6), Chapter~3), the singularity
$xy=z^3+w^3$ of Corti and Mella \cite{CM}. Each of these results has
significant applications to the Sarkisov program\index{Sarkisov!program}
and 3-fold birational geo\-metry; see paper~(6), Section~6.3 and the
forthcoming paper \cite{CM}. See also
\cite{Ka} for more examples.

 \begin{exa}
 Let $P\in Y$ be a nonsingular point and $a,b$ coprime integers; then
the\index{weighted!blowup} weighted blowup $f\colon X\to Y$ with weights
$(1,a,b)$ is an extremal divisorial\index{divisorial contraction}
contraction\index{extremal!contraction} with $f(E)=P$. Corti conjectured
in 1993 that this is the complete list.
 \end{exa}

\subsection{What is a flip?} \label{sec_flips}
 \index{flip|(}
 \begin{prb}
 Classify 3-fold flips $t\colon X\broken X^+$
 \end{prb}
In a sense, this is done in the monumental paper of Koll\'ar and Mori
\cite{KM1}, but their description is not sufficiently explicit for some
applications. The following, taken from \cite{KM1}, is the simplest
example of flips (see Brown \cite{Br} for many similar families of
examples).

 \begin{exa}
 Let $f_{m-1}(x_1, x_2)$ be a homogeneous polynomial of degree $m-1$ in 2
variables, and let $\C^\times$ act on $\C^5$ with weights
$1,1,m,-1,-1$. That is, the action is given by
 \[
 x_1, x_2, x_3, y_1, y_2 \mapsto \lambda x_1, \lambda x_2,
\la^m x_3, \la\1 y_1, \la\1 y_2.
 \]
Consider the $\C^\times$-invariant affine 4-fold $0\in A \subset \C^5$
given by
 \[
 x_4y_1=f_{m-1}(x_1,x_2).
 \]
An example of a flipping contraction\index{flipping contraction} $X\to Y$
and flip $X^+$ is obtained by taking the geometric invariant theory
quotient ($\Spec$ of the ring of invariants)
$Y=A{\mathrel{/\kern-1mm/}}\C^\times$, and setting
 \[
 X=\Proj\bigoplus_{n\le0}\Oh(nK_Y)\quad\text{and}\quad
 X^+=\Proj\bigoplus_{n\ge0}\Oh(nK_Y)
 \]
for the two sides of the flip.
 \end{exa}

Quite generally, all flips arise in this way from an affine Gorenstein
\hbox{4-fold} $0\in A=\Spec \bigoplus_{n\in \Z} \Oh(nK_Y)$ with
$\C^\times$ action (see \cite{R2}), and the problem is to write manageable
equations for $A$ in $\C^\times$-linearised coordinates. These
considerations are the starting point of Klaus Altmann's paper~(1).
 \index{flip|)}

\subsection{Mori fibre spaces}\index{Mori!fibre space}

 \begin{prb} Classify Fano 3-folds\index{Fano!3-fold} (with $B_2=1$ and
$\Q$-factorial\index{factorial@$\Q$-factorial} terminal
singularities)\index{terminal!singularities} up to biregular equivalence.
 \end{prb}

An example, and a beginning of an answer to this problem, is the list of
95 Fano\index{famous@``famous 95''} 3-fold (weighted) hypersurfaces of Reid and
Iano-Fletcher (see paper~(4), 16.6), but also the 86 codimension 2
(weighted) complete intersections (paper~(4), 16.7), the 70 codimension 3
Pfaffian cases of Alt{\i}nok \cite{Al}, Table~5.1, p.~69, etc. It is a
feature of nonsingular Fano 3-folds that they all are linear sections of
standard homogeneous varieties in their Pl\"ucker embedding. It is just
possible that, when looked at from the correct angle, this beautiful and
fundamentally simple structure will extend to singular Fano 3-folds.

\subsection{The links of the Sarkisov program} \label{sec_links}
 \index{Sarkisov!program}\index{link} \index{conic bundle|(}
It is in the study of the links of the Sarkisov program, the basis of its
applicability, that explicit birational geo\-metry really comes alive.
Here we fix a Mori fibre space $X\to S$, and we ask to classify all links
$X\broken Y$, taking off from $X$ and landing at an arbitrary Mori fibre
space $Y\to T$. For $X$ a general member of one of our famous 95
families\index{famous@``famous 95''} of Fano 3-fold hypersurfaces, this program
is carried out to completion in paper~(5). In particular, at the end of
the classification, we discover that $Y\iso X$. The links that occur fit
into a very small number of known classes. On the other hand, in doing
many concrete examples, it is our experience that each new case that we
understand involves learning how to do computations (for example, in
graded rings or in the geometry\index{projection} of projections), a
process that can sometimes be rather tricky.

\begin{exa} Consider the Fano variety $X=X_{2,2,2}\subset \PP^6$ given as
a complete intersection of three sufficiently general quadrics in $\PP^6$,
and choose a line $L\subset X$. Then there is a link $\tau_L\colon
X\broken X'$, to a conic bundle\index{conic bundle} $X'\to\PP^2$ with
discriminant curve of degree $7$.

The rational map $X\broken\PP^2$ is not too difficult to realise. For a
point $x\in \PP^6$, write $\Pi_x$ for the 2-plane spanned by $L$ and $x$. 
Write $\{Q_\la\mid \la\in \La\}$ for the net of quadrics vanishing on
$X$. Then the restriction $Q_\la{}\rest{\Pi_x}=L+\Ga_\la$ is the union of
$L$ plus a line $\Ga_\la$, and it is easy to see that $x\notin X$ if and
only if $\{\Ga_\la \mid \la\in \La\}$ is the whole of $\Pi^\vee_x$.
Mapping $x\in X$ to the quadric $Q_\la$ containing all the $\Pi_x$ (which
is unique, in general) gives a rational map $X\broken\La \iso \PP^2$,
whose fibres are conics.

The first point in seeing this process as a Sarkisov
link\index{link} is that we must understand it in explicit
biregular terms, and factor the map $X\broken\La$ as a chain of
flips,\index{flip} flops and\index{divisorial contraction} divisorial
contractions, followed by a Mori fibre space. In fact, if $x_0,\dots,x_6$
are coordinates on $\PP^6$ with $L:\{x_0=x_1=\cdots=x_4=0\}$, the equation
of $X$ can be written as
 \[
 M \begin{pmatrix} x_5 \\ x_6 \end{pmatrix}=\bq
 \]
 where $M$ is a $3\times 2$ matrix of linear forms in the variables
$x_0,\dots,x_4$ and $\bq$ is a 3-vector of quadric forms in the same
variables. Write $\pi_L\colon X\to \overline{Y} \subset \PP^4$ for the 
projection\index{projection} to $\PP^4$. The image
$\overline{Y}=\overline{Y}_4$ is the quartic 3-fold given by the equation
 \[
 \det M\bq =0.
 \]
 The singular locus of $\overline{Y}_4$ consists of the 44 ordinary nodes
 \[
 \{ \rk M\bq =1\}.
 \]
 Consequently, letting $Y\to X$ be the blow up of $L\subset X$, the
birational morphism $Y\to \overline{Y}$ contracts 44 lines with normal
bundle $(-1,-1)$, and denoting the flop $t\colon Y\broken Y'$, it is easy
to check that the rational map $Y\broken \La$ described above becomes a
morphism $Y'\to \La$, which is in fact a Mori fibre space\index{Mori!fibre
space} and a conic bundle.\index{conic bundle} This explicit construction
shows that the map $\tau_L\colon X\broken Y'$ is a Sarkisov
link\index{link} of Type~II.

Similarly, it is easy to construct self-links $\si_C\colon X\broken X$
centred on conics and twisted cubics $C\subset X$.
\end{exa}

The methods of paper~(6), Section~6.6, should prove that the only links
$X\broken Y$, starting with a general $X=X_{2,2,2}$, are the $\tau_L$ and
$\si_C$ just described. \index{conic bundle|)}

\section{Acknowledgments} This book originated in a 3-folds activity
during three months Sep--Dec~1995, organised as part of the 1995--96 EPSRC
Warwick algebraic geometry symposium. (Although it has grown almost beyond
recognition in the intervening years.) As well as the funding from EPSRC,
the 3-folds activity benefitted from a number of travel grants for Japanese
visitors from the Royal Society/ Japan Society for the Promotion of Science
administered by the Isaac Newton Institute, Cambridge and Kyoto Univ.,
RIMS, and from the support of AGE (Algebraic Geo\-metry in Europe, EU
HCM/TMR project, contract number ERBCHRXCT 940557).

It is our solemn duty, as editors, to apologise to the other authors for
our slow pace in preparing this book, which has delayed by several years
the appearance of the papers~(1)--(4). This is especially reprehensible
since the desire to understand their results and assimilate them into our
own work has been a strong motivation for our research (this applies to
Klaus Altmann's paper~(1), and more especially to Sasha Pukhlikov's
paper~(3), which provided much of the stimulus for our papers~(5) and~(6)),
and we are deeply conscious of the fact that our own work has improved
partly as a result of the delay we have inflicted on theirs. We hope that
the quality of the final product can go some way towards compensating for
our transgressions.

\bigskip
\noindent
Alessio Corti\\
DPMMS, University of Cambridge,\\
Centre for Mathematical Sciences,\\
Wilberforce Road, Cambridge CB3 0WB, U.K.\\
e-mail: a.corti@dpmms.cam.ac.uk

\medskip
\noindent
Miles Reid,\\
Math Inst., Univ. of Warwick,\\
Coventry CV4 7AL, England\\
e-mail: miles@maths.warwick.ac.uk \\
web: www.maths.warwick.ac.uk/$\sim$miles

\end{document}